\documentclass[11pt]{amsart}
\usepackage{amsfonts}
\usepackage{mathrsfs}
\usepackage{amssymb,latexsym}
\usepackage{pstcol,pstricks,color}

 \numberwithin{equation}{section}

\newtheorem{theorem}{Theorem}[section]

\newtheorem{corollary}[theorem]{Corollary}

\newtheorem{remark}[theorem]{Remark}

\newtheorem{lemma}[theorem]{Lemma}

\def\qed{\hfill $\Box$}

\title{Two kinds of hook length formulas for complete $m$-ary trees }

\begin{document}
\maketitle
\begin{center}
Yidong Sun$^\dag$ and Huajun Zhang$^\ddag$

$^\dag$Department of Mathematics, Dalian Maritime University, 116026 Dalian, P.R. China\\[5pt]
$^\ddag$Department of Mathematics, Zhejiang Normal University,
321004 Jinhua, P.R. China

{\it $^\dag$sydmath@yahoo.com.cn,\  $^\ddag$huajunzhang@zjnu.cn }
\end{center}\vskip0.5cm

\subsection*{Abstract} In this paper, we define two kinds of hook
length for internal vertices of complete $m$-ary trees, and deduce
their corresponding hook length formulas, which generalize the main
results obtained by Du and Liu.

\medskip

{\bf Keywords}: Hook length formulas, $m$-ary trees

\noindent {\sc 2000 Mathematics Subject Classification}: 05A15,
05A19

\section{Introduction}
Postnikov's hook length formula \cite{post} states that
\begin{eqnarray*}
\frac{n!}{2^n}\sum_{T}\prod_{v}\Big(1+\frac{1}{h_v}\Big)=(n+1)^{n-1},
\end{eqnarray*}
where the sum is over all unlabeled complete binary trees $T$ with
$n$ internal vertices, the product is over all internal vertices $v$
of $T$, and $h_v$ is the ``hook length" of $v$ in $T$, namely, the
number of internal vertices in the subtree of $T$ rooted at $v$.
Postnikov derived the formula indirectly and asked for a
combinatorial proof which was provided by Seo \cite{seo}, Chen and
Yang \cite{chenyang}. Later, Lascoux conjectured that
\begin{eqnarray*}
\sum_{T}\prod_{v}\Big(x+\frac{1}{h_v}\Big)=\frac{1}{(n+1)!}\prod_{i=0}^{n-1}\Big((n+1+i)x+n+1-i\Big).
\end{eqnarray*}
This is equivalent to the more suggestive form
\begin{eqnarray}\label{eqn formua1.1}
\sum_{T}\prod_{v}\frac{(h_v+1)x-h_v+1}{2h_v}=\frac{1}{n+1}\binom{(n+1)x}{n},
\end{eqnarray}
which was proved by Du and Liu \cite{duliu}. Moreover, they
generalized (\ref{eqn formua1.1}) from counting complete binary
trees to counting complete $(m+1)$-ary trees and obtained the
following formula for $(m+1)$-ary trees:
\begin{eqnarray}\label{eqn formua1.2a}
\sum_{T\in \mathcal
{T}_{n,m+1}}\prod_{v}\frac{(mh_v+1)x-h_v+1}{(m+1)h_v}=\frac{1}{mn+1}\binom{(mn+1)x}{n},
\end{eqnarray}
or equivalently
\begin{eqnarray}\label{eqn formua1.2b}
\sum_{T\in \mathcal
{T}_{n,m+1}}\prod_{v}\Big(x+\frac{1}{h_v}\Big)=\frac{x+1}{n!}\prod_{i=1}^{n-1}\Big((mn+i+1)(x+1)-(m+1)i\Big).
\end{eqnarray}
where $\mathcal {T}_{n,m+1}$ denotes the set of complete $(m+1)$-ary
trees with $n$ internal vertices, the product is over all internal
vertices $v$ of $T$.

Recall that a {\em plane forest} is a forest of plane trees that are
linearly ordered. Let $\mathcal{F}(n)$ denote the set of plane
forests with $n$ vertices. For any vertex $v$ of $F\in
\mathcal{F}(n)$, the hook length $H_v$ of $v$ is defined as the
number of vertices in the subtree rooted at $v$. Note that this
definition is slightly different to that of hook length defined
above for $(m+1)$-ary trees. Du and Liu \cite{duliu} investigated
the hook length polynomials for plane forests and obtained that
\begin{eqnarray}\label{eqn formua1.3a}
\sum_{F\in \mathcal {F}(n)}\prod_{v\in
V(F)}\Big(x+\frac{1}{H_v}\Big)=\frac{(x+1)}{n!}\prod_{i=1}^{n-1}\Big((2n+1-i)(x+1)-i\Big),
\end{eqnarray}
or equivalently,
\begin{eqnarray}\label{eqn formua1.3b}
F_n(x)=\sum_{F \in \mathcal {F}(n) }\prod_{v\in
V(F)}\frac{(2h_v-1)x-H_v+1}{H_v}=\frac{1}{2n+1}\binom{(2n+1)x}{n},
\end{eqnarray}
where $V(F)$ is the set of vertices of $F$.

It is well known that there exists a simple bijection between plane
forests and complete binary trees. For the sake of completeness, we
present it here. Given any plane forest $F\in \mathcal {F}(n)$, we
pick the first plane tree $T$ of $F$ with root $u$. Let $T'$ denote
the plane forest deduced from $T$ by removing the root $u$. Then the
bijection can be defined recursively as follows: $\psi(F)$ is the
complete binary tree with root $u$ such that it has the left subtree
$\psi(T')$ and the right subtree $\psi(F\backslash T)$.

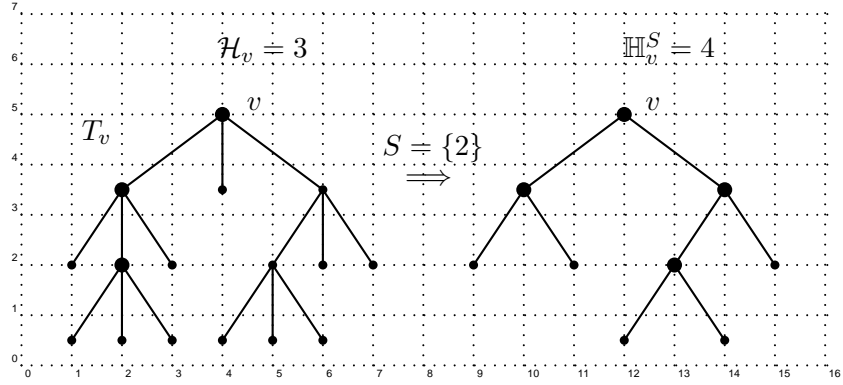
\begin{figure}[h]
\setlength{\unitlength}{0.4mm}
\begin{center}
\begin{pspicture}(11,4.5)
\psset{xunit=19pt,yunit=19pt}\psgrid[subgriddiv=1,griddots=5,
gridlabels=4pt](0,0)(16,7)

\psline(4,5)(2,3.5) \psline(4,5)(4,3.5) \psline(4,5)(6,3.5)(5,2)
\psline(2,3.5)(1,2)\psline(2,3.5)(3,2)\psline(6,3.5)(6,2)\psline(6,3.5)(7,2)
\psline(2,3.5)(2,2)(2,.5)\psline(2,2)(1,.5)\psline(2,2)(3,.5)
\psline(5,2)(4,.5)\psline(5,2)(5,.5)\psline(5,2)(6,.5)

\pscircle*(4,5){0.1}
\pscircle*(2,3.5){0.1}\pscircle*(2,2){0.1}\pscircle*(2,.5){0.06}
\pscircle*(1,.5){0.06}\pscircle*(3,.5){0.06}
\pscircle*(4,3.5){0.06}\pscircle*(6,3.5){0.06}
\pscircle*(1,2){0.06}\pscircle*(3,2){0.06}
\pscircle*(5,2){0.06}\pscircle*(6,2){0.06}
\pscircle*(4,.5){0.06}\pscircle*(5,.5){0.06}\pscircle*(6,.5){0.06}
\pscircle*(7,2){0.06}

\put(3,3.4){$v$}\put(2.6,4.1){$\mathcal {H}_v=3$}
\put(0.8,3){$T_v$}\put(5.1,2.4){$\Longrightarrow$}\put(4.8,2.8){$\tiny{S=\{2\}}$}

\psline(12,5)(10,3.5) \psline(12,5)(14,3.5)(13,2)
\psline(10,3.5)(9,2)\psline(10,3.5)(11,2)\psline(14,3.5)(15,2)

\psline(13,2)(12,.5)\psline(13,2)(14,.5)

\pscircle*(12,5){0.1} \pscircle*(10,3.5){0.1}\pscircle*(14,3.5){0.1}
\pscircle*(9,2){0.06}\pscircle*(11,2){0.06}\pscircle*(13,2){0.1}\pscircle*(15,2){0.06}
\pscircle*(14,.5){0.06}\pscircle*(12,.5){0.06}

\put(8.3,3.4){$v$}\put(8.,4.1){$\mathbb{H}_v^{S}=4$}

\end{pspicture}
\caption{The two kinds of hook length of $v$. }\label{fDD}
\end{center}
\end{figure}

It is clear that the bijection maps the hook length of $v$ in $V(F)$
to the number of internal vertices of the left component of $v$ of
$\psi(F)$. This motivates us to define the first kind of {\em hook
length} $ \mathcal {H}_v$ for an internal vertex $v$ of $m$-ary
trees $T$. Let $T_v$ denote the $m$-ary subtree of $T$ rooted at $v$
and let $T'_v$ denote the reduced tree from $T_v$ by removing the
rightmost subtree of $v$. Define $\mathcal {H}_v$ to be the number
of internal vertices of the subtree $T'_v$. See Figure \ref{fDD} for
example. Note that in the case $m=2$, the hook length $\mathcal
{H}_v$ reduces to $H_v$ up to the bijection $\psi$. Then we have the
first main result which is a generalization of (\ref{eqn
formua1.3a}) and (\ref{eqn formua1.3b}).
\begin{theorem}\label{theo 1.1}
For any integer $m\geq 2$,
\begin{eqnarray}\label{eqn 1.6}
\sum_{T\in \mathcal {T}_{n,m}}\prod_{v\in \mathcal
{I}(T)}\Big(x+\frac{1}{\mathcal
{H}_v}\Big)=\frac{(x+1)}{n!}\prod_{i=1}^{n-1}\Big((mn+1-i)(x+1)-(m-1)i\Big),
\end{eqnarray}
or equivalently,
\begin{eqnarray}\label{eqn 1.7}
\sum_{T\in \mathcal {T}_{n,m}}\prod_{v\in \mathcal
{I}(T)}\frac{(m\mathcal {H}_v-1)x-\mathcal {H}_v+1}{(m-1)\mathcal
{H}_v}=\frac{1}{mn+1}\binom{(mn+1)x}{n},
\end{eqnarray}
where $\mathcal {I}(T)$ is the set of internal vertices of $T\in
\mathcal {T}_{n,m}$.
\end{theorem}

Moreover, the definition of the first hook length inspires us
defining the second kind of hook length. Let $S$ be a subset of
$[m]=\{1,2,\dots,m\}$, for an internal vertex $v$ of $(m+1)$-ary
trees $T$, let $T_v$ denote the $(m+1)$-ary subtree of $T$ rooted at
$v$, and let $v_1,v_2,\cdots, v_{m+1}$ be the children of $v$, first
delete the subtree rooted at $v_r$ for all $r\in S$, namely delete
the $r$th subtree of $v$ for all $r\in S$; then delete the $r$th
subtree of $v_j$ for all $r\in S$ and $j\in [m+1]\setminus S$, and
then continue this process; one can obtain an $(m+1-|S|)$-ary tree
$T_v^{S}$. Define $\mathbb{H}_v^S$ to be the number of internal
vertices of $T_v^{S}$. See Figure \ref{fDD} for example. Then we
have the second main result which, in the case $S=\emptyset$,
reduces to (\ref{eqn formua1.2b}) and (\ref{eqn formua1.2a})
respectively.
\begin{theorem}\label{theo 1.2}
For any integer $m\geq 1$, $S\subset [m]$ and $s=|S|$,
\begin{eqnarray}\label{eqn formua5.1a}
\sum_{T\in \mathcal {T}_{n,m+1}}\prod_{v\in \mathcal
{I}(T)}\Big(x+\frac{1}{\mathbb
{H}_v^S}\Big)=\frac{(x+1)}{n!}\prod_{i=1}^{n-1}\Big((mn+i+1)(x+1)-(m-s+1)i\Big),
\end{eqnarray}
or equivalently,
\begin{eqnarray}\label{eqn formua5.1b}
\sum_{T\in \mathcal {T}_{n,m+1}}\prod_{v\in \mathcal
{I}(T)}\frac{\big((m-s)\mathbb {H}_v^S+1\big)x-\mathbb
{H}_v^S+1}{(m-s+1)\mathbb {H}_v^S}=\frac{1}{mn+1}\binom{(mn+1)x}{n}.
\end{eqnarray}
\end{theorem}

In the next two Sections, we present the proofs of Theorem \ref{theo
1.1} and \ref{theo 1.2} respectively.

\section{Proof of Theorem \ref{theo 1.1}}

In order to prove Theorem \ref{theo 1.1}, we need the following
lemma obtained by Seo \cite{seo}.
\begin{lemma}\label{lemma 4.2}
Fix positive integers $a$ and $b$. Let
$\Omega:=\Omega(t)=1+\sum_{n\geq 1}\Omega_nt^n$ be a formal power
series in $t$ satisfying
\begin{eqnarray*}
\Omega'=x\Omega^{b+1}+at\Omega^b\Omega',
\end{eqnarray*}
where the prime denotes the derivative of $\Omega$ with respect to
$t$. Then $\Omega_n$ can be given by
\begin{eqnarray*}
\Omega_n=\frac{x}{n!}\prod_{i=1}^{n-1}\Big(ai+bx(n-i)+x\Big).
\end{eqnarray*}
\end{lemma}
{\em Proof of Theorem \ref{theo 1.1}:} Define
\begin{eqnarray*}
\mathcal {H}_{n,m}(x)=\sum_{T\in \mathcal {T}_{n,m}}\prod_{v\in
\mathcal {I}(T)}\frac{(m\mathcal {H}_v-1)x-\mathcal
{H}_v+1}{(m-1)\mathcal {H}_v}.
\end{eqnarray*}

Given any $m$-ary tree $T\in \mathcal {T}_{n,m}$ with root $u$ for
$n\geq 1$, let $T_1,T_2,\dots,T_m$ be the $m$ subtrees of $u$ from
left to right with $i_1,i_2,\dots,i_m$ internal vertices
respectively. Then $\mathcal{H}_u=i_1+i_2+\cdots+i_{m-1}+1$.
Therefore, we can deduce the recurrence relation for $\mathcal
{H}_{n,m}(x)$,
\begin{eqnarray*}
\mathcal {H}_{n,m}(x)&=&\sum_{T\in \mathcal {T}_{n,m}}\prod_{v\in
\mathcal {I}(T)}\frac{(m\mathcal {H}_v-1)x-\mathcal
{H}_v+1}{(m-1)\mathcal {H}_v}\\
&=&\sum_{i_1+i_2+\cdots+i_{m}=n-1}\frac{(m\mathcal
{H}_u-1)x-\mathcal {H}_u+1}{(m-1)\mathcal
{H}_u}\prod_{j=1}^m\sum_{T_j\in \mathcal {T}_{{i_j},m}}\prod_{v\in
\mathcal {I}(T_j)}\frac{(m\mathcal {H}_v-1)x-\mathcal
{H}_v+1}{(m-1)\mathcal {H}_v}\\
&=&\sum_{i_1+i_2+\cdots+i_{m}=n-1}\Big(\frac{mx-1}{m-1}+\frac{1-x}{(m-1)(i_1+i_2+\cdots+i_{m-1}+1)}\Big)
\prod_{j=1}^m\mathcal {H}_{i_j,m}(x)
\end{eqnarray*}
Define the generating function for $\mathcal {H}_{n,m}(x)$ by
\begin{eqnarray*}
\mathcal {H}_{m}(x;t)=1+\sum_{n\geq 1}\mathcal {H}_{n,m}(x)t^n.
\end{eqnarray*}
Then by the above relation and the following series expansion
\begin{eqnarray*}
\mathcal {H}_{m}^k(x;t)=1+\sum_{n\geq
1}t^n\sum_{i_1+i_2+\cdots+i_{k}=n}\prod_{j=1}^k\mathcal
{H}_{i_j,m}(x),
\end{eqnarray*}
one can get
\begin{eqnarray*}
\mathcal {H}_{m}(x;t)=1+\frac{mx-1}{m-1}t\mathcal
{H}_{m}^m(x;t)+\frac{1-x}{m-1}\mathcal
{H}_{m}(x;t)\int_{0}^t\mathcal {H}_{m}^{m-1}(x;y)dy,
\end{eqnarray*}
from which, one can derive that
\begin{eqnarray*}
\mathcal {H}'_{m}(x;t)=x\mathcal {H}_{m}^{m+1}(x;t)+(mx-1)t\mathcal
{H}_{m}^m(x;t)\mathcal {H}'_{m}(x;t),
\end{eqnarray*}
where the prime denotes the derivative of $\mathcal {H}_{m}(x;t)$
with respect to $t$.

Using Lemma \ref{lemma 4.2}, we have
\begin{eqnarray*}
\mathcal
{H}_{n,m}(x)&=&\frac{x}{n!}\prod_{i=1}^{n-1}\Big((m(n-i)+1)x+(mx-1)i\Big)\\
&=&\frac{1}{mn+1}\binom{(mn+1)x}{n},
\end{eqnarray*}
which proves (\ref{eqn 1.7}).

Dividing by $\Big(\frac{1-x}{m-1}\Big)^n$ on both sides of (\ref{eqn
1.7}), and then replacing $\frac{mx-1}{1-x}$ by $x$, one can get
(\ref{eqn 1.6}) immediately. \qed\vskip0.1cm

If choose the special values $0$ or $-m$ for $x$ in (\ref{eqn 1.6}),
we get the following identities
\begin{corollary} For any integer $m\geq 2$,
\begin{eqnarray*}
\sum_{T\in \mathcal {T}_{n,m}}\prod_{v\in \mathcal
{I}(T)}\frac{1}{\mathcal {H}_v}&=&\frac{m^{n}}{mn+1}\binom{\frac{mn+1}{m}}{n},\\
\sum_{T\in \mathcal {T}_{n,m}}\prod_{v\in \mathcal
{I}(T)}\Big(m-\frac{1}{\mathcal{H}_v}\Big)&=&\frac{(m-1)^n}{n!}(mn+1)^{n-1}.
\end{eqnarray*}
\end{corollary}
\vskip0.5cm

\section{Proof of Theorem 1.2}

Define
\begin{eqnarray*}
_S\mathbb {H}_{n,m}(x)=\sum_{T\in \mathcal {T}_{n,m+1}}\prod_{v\in
\mathcal {I}(T)}\Big(x+\frac{1}{\mathbb {H}_v^S}\Big),
\end{eqnarray*}
and define the generating function for ${_S}\mathbb {H}_{n,m}(x)$ by
\begin{eqnarray*}
_S\mathbb {H}_{m}(x;t)=1+\sum_{n\geq 1} {_S}\mathbb {H}_{n,m}(x)t^n.
\end{eqnarray*}
First, we consider the case when $S$ is the empty set $\emptyset$.
Note that in this case, (\ref{eqn formua5.1a}) and (\ref{eqn
formua5.1b}) reduce to the results (\ref{eqn formua1.2b}) and
(\ref{eqn formua1.2a}) obtained by Du and Liu \cite{duliu}. Given
any $(m+1)$-ary tree $T\in \mathcal {T}_{n,m+1}$ with root $u$ for
$n\geq 1$, let $T_1,T_2,\dots,T_{m+1}$ be the $m+1$ subtrees of $u$
from left to right with $i_1,i_2,\dots,i_{m+1}$ internal vertices
respectively. Then
$\mathbb{H}_u^{\emptyset}=i_1+i_2+\cdots+i_{m+1}+1$. Therefore, we
can deduce a recurrence relation for $_\emptyset\mathbb
{H}_{n,m}(x)$,
\begin{eqnarray*}
_\emptyset\mathbb {H}_{n,m}(x)&=&\sum_{T\in \mathcal
{T}_{n,m+1}}\prod_{v\in
\mathcal {I}(T)}\Big(x+\frac{1}{\mathbb {H}_v^{\emptyset}}\Big)\\
&=&\sum_{i_1+i_2+\cdots+i_{m+1}=n-1}\Big(x+\frac{1}{\mathbb
{H}_u^{\emptyset}}\Big)\prod_{j=1}^{m+1}\sum_{T_j\in \mathcal
{T}_{{i_j},m+1}}\prod_{v\in \mathcal {I}(T_j)}\Big(x+\frac{1}{\mathbb {H}_v^{\emptyset}}\Big) \\
&=&\sum_{i_1+i_2+\cdots+i_{m+1}=n-1}\Big(x+\frac{1}{i_1+i_2+\cdots+i_{m+1}+1}\Big)
\prod_{j=1}^{m+1} {_\emptyset}\mathbb{H}_{i_j,m+1}(x).
\end{eqnarray*}

Similar to the proof of Theorem 1.1, an equation for
${_\emptyset}\mathbb {H}_{m}(x;t)$ can be derived as
\begin{eqnarray*}
{_\emptyset}\mathbb {H}_{m}(x;t)=1+xt\ {_\emptyset}\mathbb
{H}_{m}^{m+1}(x;t)+\int_{0}^t{_\emptyset}\mathbb
{H}_{m}^{m+1}(x;y)dy,
\end{eqnarray*}
from which, one can get
\begin{eqnarray}\label{eqn formua5.2}
{_\emptyset}\mathbb {H}'_{m}(x;t)=(x+1)\ {_\emptyset}\mathbb
{H}_{m}^{m+1}(x;t)+(m+1)xt\ {_\emptyset}\mathbb {H}_{m}^{m}(x;t)\
{_\emptyset}\mathbb {H}'_{m}(x;t),
\end{eqnarray}
where the prime denotes the derivative of ${_\emptyset}\mathbb
{H}_{m}(x;t)$ with respect to $t$.

For any complete $(m+1)$-ary tree $T$ with $k\geq 1$ internal
vertices and an $s$-subset $S\in [m]$, according to the definition
of the second kind of hook length, $T$ can be uniquely partitioned
into a complete $(m-s+1)$-ary tree with $n$ internal vertices for
some $n\geq 1$ and an ordered forest of $ns$ complete $(m+1)$-ary
trees. Hence we get a recurrence relation for ${_S}\mathbb
{H}_{m}(x;t)$, namely
\begin{eqnarray}\label{eqn formua5.3}
{_S}\mathbb {H}_{m}(x;t)=1+\sum_{n\geq 1}{_\emptyset}\mathbb
{H}_{n,m-s}(x)t^n {_S}\mathbb {H}_{m}^{ns}(x;t)={_\emptyset}\mathbb
{H}_{m-s}(x;{_S}\mathbb {H}_{m}^s(x;t)t).
\end{eqnarray}
Taking the derivative on both side of (\ref{eqn formua5.3}) with
respect to $t$, using (\ref{eqn formua5.2}), we have
\begin{eqnarray}\label{eqn formua5.4}
\hskip0.8cm {_S}\mathbb {H}'_{m}(x;t)=(x+1)\ {_S}\mathbb
{H}_{m}^{m+1}(x;t)+\big((m-s+1)x+s(x+1)\big)t\ {_S}\mathbb
{H}_{m}^{m}(x;t)\ {_S}\mathbb {H}'_{m}(x;t).
\end{eqnarray}
Applying Lemma \ref{lemma 4.2} to (\ref{eqn formua5.4}), one can
obtain that
\begin{eqnarray*}
{_S}\mathbb {H}_{n,m}(x)&=&\frac{x+1}{n!}\prod_{i=1}^{n-1}\Big(
\big((m-s+1)x+s(x+1)\big)i+m(x+1)(n-i)+x+1 \Big) \\
&=&\frac{(x+1)}{n!}\prod_{i=1}^{n-1}\Big((mn+i+1)(x+1)-(m-s+1)i\Big),
\end{eqnarray*}
which proves (\ref{eqn formua5.1a}).

Dividing by $\big((m-s+1)-(x+1)\big)^n$ on both sides of (\ref{eqn
formua5.1a}), and then replacing $\frac{x+1}{(m-s+1)-(x+1)}$ by $x$,
one can get (\ref{eqn formua5.1b}) immediately. \qed\vskip0.1cm

If choose the special values $m-s-1$ or $m-s$ for $x$ in (\ref{eqn
formua5.1a}), or choose $s=m$ in (\ref{eqn formua5.1b}), we get the
following identities
\begin{corollary}For any integer $m\geq 0$, $S\subset [m]$ and $s=|S|$,
\begin{eqnarray*}
\sum_{T\in \mathcal {T}_{n,m+1}}\prod_{v\in \mathcal
{I}(T)}\Big(m-s-1+\frac{1}{\mathbb{H}_v^S}\Big)&=&\frac{1}{mn+1}\binom{(mn+1)(m-s)}{n},\\
\sum_{T\in \mathcal {T}_{n,m+1}}\prod_{v\in \mathcal
{I}(T)}\Big(m-s+\frac{1}{\mathbb{H}_v^S}\Big)&=&\frac{(m-s+1)^n(mn+1)^{n-1}}{n!},\\
\sum_{T\in \mathcal {T}_{n,m+1}}\prod_{v\in \mathcal
{I}(T)}\frac{x-\mathbb {H}_v^{[m]}+1}{\mathbb
{H}_v^{[m]}}&=&\frac{1}{mn+1}\binom{(mn+1)x}{n}.
\end{eqnarray*}
\end{corollary}

\begin{remark}
Motivated by Lemma \ref{lemma 4.2} and the proof of Theorem
\ref{theo 1.2}, we can consider the function
$\Phi:=\Phi(t)=\Omega(t\Phi^s(t))$, where $\Omega(t)$ is defined in
Lemma \ref{lemma 4.2} and $\Phi=1+\sum_{n\geq 1}\Phi_nt^n$. Then it
is easy to derive that $\Phi(t)$ satisfies the following
differential equation
\begin{eqnarray*}
\Phi'=x\Phi^{b+s+1}+(a+sx)t\Phi^{b+s}\Phi',
\end{eqnarray*}
from which, by Lemma \ref{lemma 4.2}, we can deduce the explicit
expression for $\Phi_n$,
\begin{eqnarray*}
\Phi_n=\frac{x}{n!}\prod_{i=1}^{n-1}\Big(ai+bx(n-i)+(sn+1)x\Big).
\end{eqnarray*}
We wonder if there is any combinatorial explanation for the
relation.
\end{remark}
\vskip1cm

%===========================================================================
\section*{Acknowledgements} The authors are grateful to the anonymous referees for the
helpful suggestions and comments. The first author was supported by
The National Science Foundation of China.

%===========================================================================

\vskip1cm
%==============================================================================================================

\end{document}